\documentclass[a4paper]{amsart}
\usepackage{graphics}
\usepackage{color}
\usepackage{amssymb}
\usepackage[all]{xy}
\usepackage{amsmath}
\usepackage{stmaryrd}
\CompileMatrices

\newtheorem{theorem}{Theorem}[section]
\newtheorem{lemma}[theorem]{Lemma}
\newtheorem{theo}[theorem]{Theorem}
\newtheorem{proposition}[theorem]{Proposition}
\newtheorem{corollary}[theorem]{Corollary}
\theoremstyle{definition}
\newtheorem{definition}[theorem]{Definition}
\newtheorem{example}[theorem]{Example}
\theoremstyle{remark}
\newtheorem{remark}[theorem]{Remark}
\theoremstyle{remark}

\numberwithin{equation}{theorem}
\def\Chi{{\mathbb X}}

\def\Rho{{\rm R}}

\def\div{{\rm div}}

\def\KK{{\mathbb K}}

\def\ZZ{{\mathbb Z}}

\def\QQ{{\mathbb Q}}
\def\PP{{\mathbb P}}
\def\gg{\mathfrak{g}}
\def\nn{\mathfrak{n}}
\def\Of{{\mathcal{O}}}

\def\WDiv{\operatorname{WDiv}}

\def\Eff{{\rm Eff}}
\def\Ample{{\rm Ample}}
\def\SAmple{{\rm SAmple}}

\def\cov{{\rm cov}}

\def\Cl{\operatorname{Cl}}

\def\Pic{\operatorname{Pic}}
\def\Hom{{\rm Hom}}

\def\codim{{\rm codim}}

\def\Spec{{\rm Spec}}
\def\Proj{{\rm Proj}}

\def\cone{{\rm cone}}
\def\lin{{\rm lin}}

\def\cov{{\rm cov}}
\def\SL{{\rm SL}}

\def\Ker{{\rm Ker}\,}

\newcounter{itemnumber}

\begin{document}
\sloppy
\author[Ivan V.~Arzhantsev]%
{Ivan V. Arzhantsev} 
\title[Projective embeddings with small boundary]
{Projective embeddings of homogeneous spaces\\ 
with small boundary}
\thanks{Supported by INTAS YS 05-109-4958}
\address{Department of Higher Algebra, Faculty of Mechanics and Mathematics, Moscow State University,
Leninskie Gory, GSP-2 Moscow, 119992, Russia}
\email{arjantse@mccme.ru}
\begin{abstract}
We study open equivariant  projective embeddings of homogeneous spaces such that the complement
of the open orbit has codimension $\ge 2$. Criterions of existence of such an embedding 
are considered and finiteness of isomorphism classes of such embeddings for a given homogeneous 
space is proved. Any embedding with small boundary is realized as a GIT-quotient associated with
a linearization of the trivial line bundle on the space of the canonical embedding. The generalized Cox's
construction and the theory of bunched rings allow us to describe basic geometric properties of
embeddings with small boundary in combinatorial terms. 
\end{abstract}

\maketitle

%%%%%%%%%%%%%%%%%%%%%%%%%%%%%%%%%%%%%%%%%%%%%%%%%%%%%%%%%%%%%%%%%%%%%%%%

\section*{Introduction}
 Let $G$ be a connected affine algebraic group over an algebraically closed field $\KK$ of
characteristic zero, and $H$ be a closed subgroup of $G$. 
The main object of the paper are equivariant embeddings of the homogeneous space $G/H$, i.e.,
open dense $G$-equivariant embeddings $i:G/H\hookrightarrow X$, where $X$ is a normal
$G$-variety. The theory of embeddings of homogeneous spaces is a well-developed branch of algebraic
transformation group theory. 
A major part of classification results here is based on the approach proposed in~\cite{LuVu} 
and known as the Luna-Vust theory. In principle, this theory describes all embeddings of a
given homogeneous space $G/H$, but such a description is constructive only for spaces with 
complexity $\le 1$.

 Bearing this restriction in mind, one may try to classify embeddings and to study their properties 
under some conditions on the 
variety $X$. For example, if $X$ is affine, we get additional technical tools: $G$-module structure
on the algebra of regular functions on $X$ and interaction of this structure with multiplication.
A survey of recent results on affine embeddings is given in~\cite{Ar1}. 
Another possible restriction is to assume that $X$ is projective. But it turns out that this case is
in no sense simpler than the general one. In this work we study a special class of projective
embeddings. Namely,  {\it a projective embedding with small boundary} is an embedding
$i:G/H\hookrightarrow X$ such that $X$ is normal projective and the boundary 
$X\setminus i(G/H)$ does not contain divisors. Such embeddings appeared in earlier papers
(see \cite{BBII}, \cite{BBK}, \cite[Section 23 B]{Gr}), but, as far as we know, their first systematic
investigation was undertaken in~\cite{ArHa}.  The present paper is a direct continuation of this investigation.  
Note that in \cite{ArHa} we gave a combinatorial description of a wide class of so-called $A_2$-maximal
embeddings with small boundary, while here we deal with projective embeddings only.

 It is well known that the algebra of regular functions $\Of(X)$ on an irreducible projective variety $X$ 
consists of constants. Thus the existence of a projective $G/H$-embedding with small boundary
implies $\Of(G/H)=\KK$. Recall that a closed subgroup $H$ of an algebraic group $G$ is said to be {\it epimorphic}
if $\Of(G/H)=\KK$. Some characterizations, interesting properties and examples of epimorphic subgroups may be found in 
\cite{Bi}, \cite{BBI}, \cite{BBII}, \cite{BBK}, \cite[Section 23 B]{Gr}. 

 In Section~\ref{sec:1}, it is shown that a projective embedding of $G/H$ 
with small boundary is determined by a character of $H$. Using a subgroup from~\cite{BBII} (see also \cite{BBK}) 
and Nagata-Steinberg's counterexample to Hilbert's 14th Problem,  
we give an example of homogeneous space $G/H$ with epimorphic $H$ that does not admit 
completions with small boundary (Section~\ref{sec:2}). This example shows that the condition "$H$ is epimorphic"
is not sufficient for the existence of complete embedding with small boundary.

A closed subgroup $H$ of $G$ is called {\it observable} if the homogeneous space
$G/H$ is quasiaffine. Further, a subgroup $H\subseteq G$ is a {\it Grosshans} subgroup if
$H$ is observable and the algebra $\Of(G/H)$ is finitely generated. It is known that 
a homogeneous space admits an affine embedding with small boundary $G/H\hookrightarrow Z$ 
if and only if $H$ is a Grosshans subgroup \cite{Gr}. In this case a 
(normal) affine embedding of $G/H$ with small boundary is unique (up to $G$-equivariant
isomorphism), and the variety $Z$ is the spectrum of the algebra $\Of(G/H)$.
The embedding $G/H\hookrightarrow Z=\Spec(\Of(G/H))$ is called the {\it canonical} embedding of $G/H$.

It is natural to ask, can these results be transfered to the projective case. 
Here an analogue of a Grosshans subgroup may be defined as an epimorphic subgroup obtained from a
Grosshans subgroup by a torus extension (in~\cite{ArHa} we call it a {\it Grosshans extension}). 
Corresponding homogeneous spaces admit projective embeddings with small boundary. Concerning uniqueness,
we show that the number 
of projective $G/H$-embeddings with small boundary,
where $H$ is a Grosshans extension, is finite (Section~\ref{sec:3}). Moreover, any projective embedding with small boundary
may be realized as a categorical quotient of the set of semistable points of a linearized trivial line bundle
over the canonical embedding $G/H_1\hookrightarrow Z$ with respect to a torus action,
where $H_1$ is the intersection of kernels of all characters of $H$. In contrast to~\cite{ArHa}, here we deal only with 
elementary facts of Geometric Invariant Theory (GIT). This allows us to take off conditions on $G$ and $H$.  

Using the notion of the total coordinate ring and a generalization of Cox's construction from toric geometry,
it was shown in~\cite{ArHa} that under some conditions on the pair $(G,H)$ projective embeddings of $G/H$
with small boundary are parametrized by "interior" cones of a fan $\Sigma(G/H)$ that appears as the GIT-fan 
of an action of a torus on an affine factorial variety. Moreover, equivariant morphisms 
between embeddings correspond to the face relation on the set of cones of $\Sigma(G/H)$. 
These results are discussed in Section~\ref{sec:30}.

 The theory of bunched rings developed in~\cite{BH1} provides a combinatorial description of basic geometric
properties of normal varieties with a free finitely generated divisor class group and a finitely generated
Cox ring (compare~\cite{HK}). In Section~\ref{sec:31}, we reformulate these results for projective embeddings with
small boundary and describe the Picard group, the cones of effective, semiample and ample divisors, characterize
locally factorial and $\QQ$-factorial embeddings. If the total space of the quotient morphism is smooth, 
smoothness of the embedding turns out to be equivalent to local factoriality. If the space $Z$
is an "intrinsic complete intersection", then the canonical class of the embedding may be calculated effectively. 

 In the last section we deal with examples. In particular, projective embeddings with small boundary for 
$G=\SL(3)$ are described and an epimorphic subgroup $H$ of maximal rank in  
$G=\SL(4)$ such that $G/H$ admits many projective embeddings with small boundary is given.  
In this case we compute explicitly the graph of equivariant morphisms.

 The author is indebted to J.~Hausen for fruitful collaboration and 
numerous useful discussions. 

%%%%%%%%%%%%%%%%%%%%%%%%%%%%%%%%%%%%%%%%%%%%%%%%%%%%%%%%%%%%%%%%%%%%%%%

\section{Epimorphic subgroups and characters}
\label{sec:1}

The following result provides a criterion of existence of a projective embedding with small boundary for a homogeneous
space $G/H$. It is known to specialists (some its variant may be found in~\cite[Thm.~1]{BBII}), but for convenience
of the reader we give here a complete proof. 

\begin{theo} \label{T1}
Let $G$ be a connected affine algebraic group and $H$ be a closed subgroup of $G$.
The following conditions are equivalent:
\begin{enumerate}
\item[(1)] there exists a projective embedding $G/H\hookrightarrow X$ with small boundary;
\item[(2)] the subgroup $H$ is epimorphic and there is a character $\chi$ of $H$ such that $\Ker(\chi)$ is
a Grosshans subgroup of $G$.
\end{enumerate}
\end{theo}

\begin{proof}
Let $G/H\hookrightarrow X$ be a projective embedding with small boundary. 
Then $H$ is epimorphic and there exists a finite-dimensional rational $G$-module $V$ and a 
closed equivariant embedding $X\subseteq\PP(V)$ \cite[Th.~1.7]{PV89}. 
Let $\tilde{X}\subseteq V$ be the affine cone over the image of $X$. Taking a composition of the embedding
$X\subseteq\PP(V)$ with a power of Veronese embedding, one may assume that
the cone $\tilde{X}$ is normal.

\begin{lemma} \label{L1}
The action $G:\tilde{X}$ has an open orbit.
\end{lemma}

\begin{proof}
Consider $x_0\in X$ with the stabilizer $G_{x_0}$ coinciding with $H$. Take a non-zero
vector $\tilde{x_0}$ on the line $x_0$ in $V$. If its orbit
is not open in $\tilde{X}$, then the stabilizer $\tilde{H}$ of the point
$\tilde{x_0}$ has finite index in $H$. The homogeneous space
$G/\tilde{H}\cong G\tilde{x_0}$ is quasiaffine. Hence $G/H$ is quasiaffine~\cite[Cor.2.2]{Gr}.
But $H$ is epimorphic, a contradiction.
\end{proof}

Let $\chi$ be a character of $H$ such that $h\cdot\tilde{x_0}=\chi(h)\tilde{x_0}$ for
any $h\in H$. Then the stabilizer $G_{\tilde{x_0}}=\tilde{H}$ coincides with $\Ker(\chi)$.
Moreover, the orbit $G\tilde{x_0}$ in $\tilde{X}$ is conic and the embedding 
$G/\tilde{H}\hookrightarrow\tilde{X}$, $g\tilde{H}\to g\cdot\tilde{x_0}$ is a 
(normal) affine embedding of $G/\tilde{H}$ with small boundary. This implies that the algebra
$\Of(G/\tilde{H})=\Of(\tilde{X})$ is finitely generated and thus
$\tilde{H}=\Ker(\chi)$ is a Grosshans subgroup of $G$.

\smallskip

Conversely, let $\chi:H\to\KK^{\times}$ be a character with the kernel being a Grosshans
subgroup of $G$. Since $H$ is epimorphic, the quotient $F:=H/\Ker(\chi)$ is a one-dimensional
torus, and its $G$-equivariant action on $G/\Ker(\chi)$ by right translation
defines a positive grading on the algebra $\Of(G/\Ker(\chi))$. (Recall that a
$\ZZ$-grading $A=\oplus_{n\in\ZZ} A_n$ on a $\KK$-algebra $A$ is positive
if $A_0=\KK$ and $A_n=0$ for any $n<0$.) Consider a projective $G$-variety 
$X:=\Proj(\Of(G/\Ker(\chi)))$ defined by this grading. Let 
$$
G/\Ker(\chi)\hookrightarrow Z:=\Spec(G/\Ker(\chi))
$$ 
be the canonical embedding of $G/\Ker(\chi)$, and $O$ be the $G$-fixed point on $Z$
corresponding to the maximal ideal of positive components $\oplus_{n>0}\Of(G/\Ker(\chi))_n$
in $\Of(G/\Ker(\chi))$. There is a canonical $G$-equivariant surjection 
$p:Z\setminus\{O\}\to X$, whose fibers are $F$-orbits. The image of
the open orbit in $Z$ is an open orbit in $X$ isomorphic to $G/H$. Since
all $F$-orbits on $Z\setminus\{O\}$ are one-dimensional and the boundary of $X$
is the image of the boundary of $Z$, the embedding $G/H\hookrightarrow X$ is a (normal)
projective embedding with small boundary.
\end{proof}

\begin{remark}
Condition (2) of Theorem~\ref{T1} is mentioned in~\cite[Section 23 B]{Gr} as property (FG)
of a subgroup $H$.  
\end{remark}

\begin{remark}
We assume $G$ to be connected. However this restriction is not essential:
projective embeddings $G/H\hookrightarrow X$  with small boundary
are in bijection with corresponding embeddings of $G^0/(G^0\cap H)$. Indeed, normality of $X$ 
implies that different irreducible components of $X$ do not intersect.  
\end{remark}

\begin{corollary}[of the proof] \label{C1}
Any projective embedding $G/H\hookrightarrow X$ with small boundary may be obtained as
$X=\Proj(\Of(G/\Ker(\chi)))$ for some character $\chi$ of the subgroup $H$. 
\end{corollary}

Let us denote the embedding corresponding to a character $\chi$ as $G/H\hookrightarrow X(\chi)$,
or just $X(\chi)$. In this context, two natural questions arise:

\smallskip

\ (Q1) \ May one describe "constructively" the characters of $H$ which 
define projective embeddings with small boundary ?

\smallskip

\ (Q2) \ When do two characters $\chi_1$ and $\chi_2$ define $G$-isomorphic embeddings ?

\smallskip

To answer the first question, we start with a criterion of observability of the kernel $\Ker(\chi)$ in $G$.
Let us consider a more general problem.
Let $H$ be a closed subgroup of $G$ and $\Chi(H)$ be the group of characters of $H$. Set
$$
H_1:=\cap_{\chi\in\Chi(H)}\Ker(\chi). 
$$
This subgroup is observable in $G$ \cite[Prop.~3.13]{ArHa}, the factor group 
$H/H_1$ is diagonalizable  and thus is isomorphic to the direct product of a torus $Т$ and a finite abelian group $A$.
Consider an intermediate subgroup $H_1\subseteq H'\subseteq H$. We are going to determine when $H'$
is observable or epimorphic in $G$. Since these properties depend only on the connected component of $H'$, 
we may assume that the image $\phi(H')$ under the projection $\phi:H\to H/H_1$ is a subtorus 
$S$ in the torus $T$. Such a subtorus is defined by a primitive sublattice $\Rho(S)\subseteq\Chi(T)$,
$\Rho(S):=\{\chi\in\Chi(T) : \chi(s)=1 \ \forall s\in S\}$.   

 The torus $T$ acts $G$-equivariantly on the homogeneous space $G/H_1$ 
be right translation. This action defines a $G$-invariant grading

\begin{equation}
\label{one}
\Of(G/H_1)=\bigoplus_{\mu\in\Chi(T)}\Of(G/H_1)_{\mu}, 
\end{equation}

\noindent где $\Of(G/H_1)_{\mu}:=\{f\in\Of(G/H_1) : f(ghH_1)=\mu(\phi(h))f(gH_1) \ \forall \ g\in G, h\in\phi^{-1}(T) \}.$

\smallskip

Consider a semigroup $\Chi(G/H_1, T):=\{\mu\in\Chi(T) : \Of(G/H_1)_{\mu}\ne 0\}$ and a cone $C=C(G/H_1, T)$ that
is the closure of the cone generated by $\Chi(G/H_1,T)$ in the space $\Chi(T)_{\QQ}:=\Chi(T)\otimes_{\ZZ}\QQ$. 
(If the semigroup $\Chi(G/H_1, T)$ is finitely generated, it generates a closed polyhedral cone.)
Since the space $G/H_1$ is quasiaffine, the $T$-action on $\Of(G/H_1)$ is effective and the cone $C$ has maximal dimension
in $\Chi(T)_{\QQ}$. Let $C^{\circ}$ be the interior of $C$.

\begin{proposition}
Let $H_1\subseteq H'\subseteq H$ and $S=\phi(H')$. Then
\begin{enumerate}
\item$$
\Of(G/H')=\bigoplus_{\mu\in\Chi(G/H_1, T)\cap\Rho(S)} \Of(G/H_1)_{\mu};
$$
\item if $C^{\circ}\cap\Rho(S)\ne\emptyset$, then the subgroup $H'$ is observable in $G$; 
if $H_1$ is a Grosshans subgroup of $G$, then the converse is true;
\item the subgroup $H'$ is epimorphic in $G$ if and only if $H$ is epimorphic and $\Rho(S)\cap\Chi(G/H_1,T)=\{0\}$.
\end{enumerate}
\end{proposition}

\begin{proof}
Statement (i) stems from formula~\ref{one} and the equality $\Of(G/H')=\Of(G/H_1)^S$.

To verify (ii), let $D\subseteq\Of(G/H')$ be a finitely generated $(G\times T/S)$-invariant
subalgebra in $\Of(G/H')$ with the quotient field $QD$ coinciding with $Q\Of(G/H')$. One has an affine
embedding $G/H''\hookrightarrow\Spec(D)$, where $H''$ is an observable subgroup of $G$ 
containing $H'$. Moreover, $H'$ is observable if and only if $H'=H''$.

Let $B\subseteq\Of(G/H_1)$ be a finitely generated $(G\times T)$-invariant subalgebra with
$QB=Q\Of(G/H_1)$. One has an affine embedding $G/H_1\hookrightarrow\Spec(B)$.
Extending $B$, one may assume that $D\subset B$. 
Since the algebra of invariants $B^S$ is finitely generated, we may replace $D$ by $B^S$.
Then the embedding $D\subset B$ corresponds to the $S$-quotient morphism $p:\Spec(B)\to\Spec(D)$. 
All fibers of the first arrow in $G/H_1\to G/H'\to G/H''$ are isomorphic to $S$,
and in order to prove that $H'=H''$ it is sufficient to show that the general fiber of the morphism $p$ is an 
$S$-orbit, or, equivalently, the cone generated by the weight semigroup $\Chi(\Spec(B), S)$
coincides with the space $\Chi(S)_{\QQ}$. The semigroup $\Chi(\Spec(B), S)$ is the image
of $\Chi(\Spec(B), T)$ under the restriction $\Chi(T)\to\Chi(S)$ of characters to the subtorus.
The condition $C^{\circ}\cap\Rho(S)\ne\emptyset$ means that the kernel $\Rho(S)$ of this projection 
intersects the interior of the cone $C$, and thus $C$ projects to $\Chi(S)_{\QQ}$ surjectively.

Conversely, if $H_1$ is a Grosshans subgroup, one may put $B=\Of(G/H_1)$, and
$G/H_1\hookrightarrow\Spec(B)$ is the canonical embedding of $G/H_1$. The condition $H'=H''$
shows that the general fiber of the morphism $p$ contains $S$ as an open orbit. If
this orbit is closed, we come again to the conclusion that the projection
$C\to\Chi(S)_{\QQ}$ is surjective, or, equivalently, $C^{\circ}\cap\Rho(S)\ne\emptyset$.
If an open $S$-orbit in general fiber of $p$ is not closed, then its closure contains an 
$S$-orbit of dimension $\dim S-1$. $G$-translations of this orbit form a $G$-orbit $O$ in $\Spec(B)$,
which is mapped by $p$ to $G/H''$ surjectively. This implies that
$\dim O=\dim G/H_1-1$, a contradiction with $\codim_{\Spec(B)}(\Spec(B)\setminus (G/H_1))\ge 2$.

For (iii), note that $H'$ is epimorphic in $G$ if and only if
$\Of(G/H_1)_0=\KK$, i.e., $H$ is ephimorphic, and $\Of(G/H')=\Of(G/H_1)_0$. The last condition is equivalent to
$\Rho(S)\cap\Chi(G/H_1,T)=\{0\}$.
\end{proof}

Any character $\mu\in\Chi(T)$ may be extended to a character of $H/H_1$ (by setting it equals $1$ on $A$),
and thus to a character $\chi\in\Chi(H)$. Conversely, any character $\chi\in\Chi(H)$ is trivial on $H_1$,
hence defines a character of $H/H_1$ and a character $\mu\in\Chi(T)$. Since the connected components of the kernels of 
proportional characters coincide, we get

\begin{corollary}
Let $H$ be an epimorphic subgroup of $G$, $\chi\in\Chi(H)$, and $\mu\in\Chi(T)$ is the character corresponding to
$\chi$. Then 
\begin{enumerate}
\item if $\mu\in C^{\circ}$, then the subgroup $\Ker(\chi)$ is observable in $G$; if $H_1$
is a Grosshans subgroup, then the converse is true;
\item $\Ker(\chi)$ is epimorphic in $G$ if and only if  
$m\mu\notin\Chi(G/H_1,T)$ for any $m>0$.  
\end{enumerate}
\end{corollary}

\begin{example}
Let $G$ be a connected reductive group and $H=B=TB^u$ be a Borel subgroup of $G$ with a maximal torus 
$T$ and the maximal unipotent subgroup $B^u$. Then $H_1=B^u$, 
$\Chi(G/H_1,T)$ is the semigroup of dominant weights, and
$C$ coincides with the positive Weyl chamber. In this case, for $\chi\in\Chi(B)$ the subgroup $\Ker(\chi)$
is observable in $G$ if and only if the weight $\chi$ is strictly dominant, and $\Ker(\chi)$ is epimorphic if and only
if $\chi$ is not dominant.
\end{example}
 
 Suppose now that the subgroup $\Ker(\chi)$ is observable. The question when  
$\Ker(\chi)$ is a Grosshans subgroup seems to be very difficult. A classification of
Grosshans subgroups is connected directly with Hilbert's 14th Problem and is very far 
from being complete. In particular, we do not know answers to the following questions:

\smallskip

\ (Q3) \ Consider characters $\chi_1,\chi_2\in\Chi(H)$ with $\Ker(\chi_1)$ and
$\Ker(\chi_2)$ being observable in $G$. May it turn out that $\Ker(\chi_1)$ is a Grosshans
subgroup, but $\Ker(\chi_2)$ is not ? 

\smallskip

\ (Q4) \ Suppose that for some $\chi\in\Chi(H)$ the kernel $\Ker(\chi)$ is a Grosshans 
subgroup of $G$. Does it follow that $H_1$ is a Grosshans subgroup of $G$ ?

\smallskip

(Remark that the positive answer to (Q4) implies the negative answer to (Q3).)

\smallskip

Nowadays some sufficient conditions for a subgroup to be Grosshans and some examples of
observable non-Grosshans subgroups are known, see a survey of these results in~\cite{Gr}. Using one of the examples and following~\cite{BBII}, 
we give in the next section an example of an epimorphic subgroup $H$ in a semisimple
group $G$ such that for any $\chi\in\Chi(H)$ the subgroup $\Ker(\chi)$ in not Grosshans.

\smallskip

We finish this section with a condition on stabilizers of points on projective embeddings with small boundary. 

\begin{proposition}
Let $G$ be a connected reductive group and $G/H\hookrightarrow X$ be a projective embedding with small boundary.
Then for any point $x\in X$ the stabilizer $G_x$ is contained in a proper parabolic subgroup of $G$.
\end{proposition}

\begin{proof}
If $G_x$ is not contained in a proper parabolic subgroup, then the point $x$ has a $G$-invariant
affine neighborhood $U$ in $X$ \cite[Th.4.1]{Tr}. But then $G/H\subset U$ and the space
$G/H$ quasiaffine, a contradiction. 
\end{proof}

\begin{corollary}
Let $G$ be a connected reductive group and $G/H\hookrightarrow X$ be a projective embedding with small boundary.
Then $X$ does not contain $G$-fixed points. 
\end{corollary}

Remark that for non-reductive $G$ the last statement is not true: 
one may consider the action of a maximal parabolic subgroup 
$P\subset\SL(3)$ on $\PP^2$. 

%%%%%%%%%%%%%%%%%%%%%%%%%%%%%%%%%%%%%%%%%%%%%%%%%%%%%%%%%%%%%%%%%%%%%%%%

\section{Nagata-Steinberg's Counterexample}
\label{sec:2}

Put $G=\SL(2)\times\dots\times\SL(2)$ (9 copies)
and fix numbers $a_1,\dots,a_9$ such that $\sum_{i=1}^9 a_i\ne 0$
and the points $(a_1^2,a_1^3), (a_2^2,a_2^3), (a_3^2,a_3^3)$ do not lie on a line. 
Define a subgroup $H\subset G$ as
$$
H=\left\{
\begin{pmatrix}
t & c_1 \\
0 & t^{-1} 
\end{pmatrix},\dots,
\begin{pmatrix}
t & c_9 \\
0 & t^{-1} 
\end{pmatrix}
\right\},
$$
\noindent where $t\in\KK^{\times}$, $c_1,\dots,c_9\in\KK$, and $\sum_{i=1}^9 c_i=0$,
$\sum_{i=1}^9 a_i^2c_i=0$, $\sum_{i=1}^9 a_i^3c_i=0$. 

\smallskip

This subgroup is a semidirect product of a one-dimensional torus and a six-dimensional commutative unipotent group.

\begin{lemma} \label{L2}
The subgroup $H$ is epimorphic in $G$.
\end{lemma}

\begin{proof}
Assume the converse. Then $H$ is contained in a proper subgroup $F$
that is observable in $G$ \cite[Lemma~23.5]{Gr}. The subgroup $F\subset G$ is observable if and only if
either it is reductive or  it stabilizes a highest weight vector in some non-trivial 
simple $G$-module \cite[Lemma~7.7]{Gr}. A unique
(up to conjugation) proper connected reductive subgroup of
$\SL(2)$ is a maximal torus. Considering subgroups that project to each copy of 
$\SL(2)$ surjectively, one shows that $H$ is not contained in a proper
reductive subgroup of $G$. Therefore $H$ should stabilizes a highest weight vector. 

Any simple $G$-module is isomorphic to $V=V_1\otimes\dots\otimes V_9$, where 
$V_i$ are simple $\SL(2)$-modules of dimension $d_i$, and $G$ acts on $V$ component-wise.
Let $v=v_1\otimes\dots\otimes v_9$ be a highest weight vector stabilized by $H$. The lines $\langle v_i\rangle$
are preserved by the standard Borel subgroup of $\SL(2)$, thus $v_i$ is multiplies by
$t^{d_i-1}$, and $v$ is multiplied by $t^{d_1+\dots+d_9-9}$. But $d_i\ge 1$ and
at least one $d_i>1$, a contradiction with $H\subseteq G_v$.
\end{proof}

\begin{lemma} \label{L3}
For any $\chi\in\Chi(H)$ the kernel $\Ker(\chi)$ is not a Grosshans subgroup of $G$.
\end{lemma}

\begin{proof}
Consider a character $\chi_1\in\Chi(H)$, $\chi_1(h):=t$. The kernel $\Ker(\chi_1)$ coincides with the unipotent radical $H^u$
of the subgroup $H$. If $H^u$ is a Grosshans subgroup, then for any $G$-module $V$ of algebra of $H^u$-invariants
$\Of(V)^{H^u}$ is finitely generated \cite[Th.9.3]{Gr}. However if one considers the $G$-module $V=\KK^2\oplus\dots\oplus\KK^2$ (9 copies)
with the component-wise $G$-action, then the algebra $\Of(V)^{H^u}$ is not finitely generated \cite{St}.

Now take an arbitrary character $\chi_n(h):=t^n$. The algebra $\Of(G/H^u)$ is the integral closure of the algebra
$\Of(G/\Ker(\chi_n))$ in the field $\KK(G/H^u)$, which is a finite extension of the field $\KK(G/\Ker(\chi_n))$.
Hence if  $\Of(G/\Ker(\chi_n))$ is finitely generated, then $\Of(G/H^u)$ is finitely generated, a contradiction.
\end{proof}

\begin{theo} \label{T2}
The homogeneous space $G/H$ admits no embeddings $G/H\hookrightarrow X$ with small boundary,
where the variety $X$ is complete.
\end{theo}

\begin{proof}
Theorem~\ref{T1} and Lemma~\ref{L3} imply that $G/H$ admits no projective embeddings with small boundary.

Let us come to an arbitrary completion $G/H\hookrightarrow X$ with small boundary. Note that $\Of(X)=\Of(G/H)=\KK$
(Lemma~\ref{L2}) and $\Cl(X)=\Cl(G/H)$, where $\Cl(X)$ is the divisor class group of a variety $X$. It is known
that if $G$ is connected simply connected semisimple and $H$ is a closed subgroup of $G$,
then $\Cl(G/H)\cong\Chi(H)$ \cite[Th.~4]{Po}. In our case, one has $\Cl(X)\cong\ZZ$.

\begin{proposition}
Let $X$ be a normal variety with $\Of(X)=\KK$ and $\Cl(X)\cong\ZZ$. Then
$X$ is quasiprojective.
\end{proposition}

\begin{proof}
Let $X=\cup_i U_i$ be an affine covering of $X$. The complement $X\setminus U_i$ is a union of 
finitely many prime divisors $D_{i1}\cup\dots\cup D_{ik_i}$. Consider a divisor $D$ whose class
$[D]$ generates $\Cl(X)$. Then $[D(i)]=s_i[D]$ for some integers
$s_i$, where $D(i):=D_{i1}+\dots+D_{ik_i}$. All integers $s_i$ have that same sign. 
Indeed, if, for example, $s_1\le 0$, $s_2\ge 0$, then $s_2D(1)-s_1D(2)$ is a principal effective divisor.
The condition $\Of(X)=\KK$ implies $s_2D(1)-s_1D(2)=0$, hence $X\setminus U_1=\emptyset$ and $X$ is affine, a contradiction.

We may assume that all $s_i>0$. Replacing the divisors $D(i)$ by their multiples, we may assume that they are
linearly equivalent. Since the complements to the supports of all $D(i)$ form an affine covering of $X$,
any $D(i)$ is ample.
\end{proof}

So, if $G/H\hookrightarrow X$ is a completion with small boundary, then the variety $X$ is quasiprojective and complete.
Thus $X$ is projective, but this is impossible. Theorem~\ref{T2} is proved. 
\end{proof}

%%%%%%%%%%%%%%%%%%%%%%%%%%%%%%%%%%%%%%%%%%%%%%%%%%%%%%%%%%%%%%%%%%%%%%%%%%%%

\section{Classification of projective embeddings with small boundary}
\label{sec:3}

In this section we give a combinatorial classification of projective $G/H$-embeddings with small boundary
under the assumption that the subgroup $H_1$ is a Grosshans subgroup of $G$. 
Since for any character $\chi\in\Chi(H)$ the algebra $\Of(G/\Ker(\chi))$ is a subalgebra of
$\Of(G/H_1)$ consisting of functions semiinvariant with respect to the quasitorus $\Ker(\chi)/H_1$, 
finite generation of $\Of(G/H_1)$ implies finite generation of $\Of(G/\Ker(\chi))$. This shows that
a character $\chi$ defines a projective $G/H$-embedding with small boundary if and only if
the subgroup $\Ker(\chi)$ is observable in $G$.

 Consider the affine $G$-variety $Z:=\Spec(\Of(G/H_1))$. Recall that the factor group $H/H_1$ is isomorphic to a
direct product $T\times A$, where $T$ is a torus and $A$ is a finite abelian group. We are interested in characters 
$\chi\in\Chi(H)$ that define projective $G/H$-embeddings with small boundary, and characters $\chi$ and $n\chi$ define
isomorphic embeddings. Hence we may consider characters  which are trivial on $A$ and identify them with characters 
of $T$. The torus $T$ acts $G$-equivariantly on $G/H_1$, and thus on
$\Of(G/H_1)$ and on $Z$. Let $f_1,\dots,f_m$ be a generating system of the algebra $\Of(Z)$ consisting of $T$-semiinvariants,
i.e., $t\cdot f_i=\mu_i(t)f_i$ for some $\mu_1,\dots,\mu_m\in\Chi(T)$. Consider a semigroup 
$\Chi(Z,T)=\Chi(G/H_1,T)$ consisting of weights $\mu$ such that the homogeneous component $\Of(Z)_{\mu}$ is non-zero.
Let $C=C(G/H_1, T)$ be the cone generated by this semigroup. Clearly, the semigroup $\Chi(Z,T)$ is generated by the weights
$\mu_1,\dots,\mu_m$, and $C$ is a convex polyhedral cone of maximal dimension in $\Chi(T)_{\QQ}$. 

 Our aim is to realize projective embeddings of $G/H$ with small boundary as GIT-quotients corresponding to various $T$-linearizations
of the trivial line bundle on the affine variety $Z$. Here we use some results from \cite[Section~2]{BH}.
With any point $z\in Z$ one associates an {\it orbit semigroup}, i.e. the semigroup of weights $\mu\in\Chi(Z,T)$ such that
there is a semiinvariant $f\in\Of(Z)_{\mu}$ with $f(z)\ne 0$. Define an {orbit cone} of a point $z$ as the cone
$\omega(z)\subseteq\Chi(T)_{\QQ}$ generated by its orbit semigroup. One may check that the cone $\omega(z)$ 
is generated by $\mu_i$ with $f_i(z)\ne 0$. In particular,  the collection of orbit cones is finite. Further, with any character
$\chi\in C\cap\Chi(T)$ one associates a {\it GIT-cone} $\sigma(\chi)$:
$$
\sigma(\chi):=\bigcap_{z\in Z, \chi\in\omega(z)} \omega(z). 
$$

Recall that a finite collection $\Sigma$ of convex polyhedral cones in a finite dimensional rational vector space конечномерном $V$ 
is said to be a {\it fan}, if  

(i) all faces of an element of $\Sigma$ belong to $\Sigma$; 

(ii) the intersection of any two elements of $\Sigma$ is a face of each of them. 

\smallskip

(Sometimes all cones in a fan are supposed to be strictly convex. Here we omit this condition. In the situation studied below
it is fulfilled automatically.)
The {\it support} of a fan $\Sigma$ is the set of vectors $v\in V$ which are contained in a cone from $\Sigma$.
 
\begin{proposition} \cite[Th.~2.11]{BH} \label{P1}
The set of GIT-cones
$$
\Sigma(Z):=\{\sigma(\chi) : \chi\in C\cap\Chi(T)\}
$$ 
is a fan with the support $C$. 
\end{proposition}

The fan $\Sigma(Z)$ is called the {\it GIT-fan} of an affine $T$-variety $Z$. This fan may be calculated effectively.
For example, if the algebra $\Of(Z)$ is given in terms of generators and relations, an algorithm 
that computes $\Sigma(Z)$ may be found in \cite[Remark~1.3]{ArHa1}.

Any character $\chi\in\Chi(T)$ defines a $T$-linearization of the trivial line bundle on $Z$:
$$
T\times (Z\times\KK)\to Z\times \KK, \ \ (t, z, a)\to (t\cdot z, \chi(t)a),
$$
and the set of semistable points of this linearization is
$$
Z^{ss}(\chi):=\{z\in Z : f(z)\ne 0 \ \text{for} \ \text{some} \ f\in\Of(Z)_{s\chi}, s>0\}. 
$$
Clearly, the subset $Z^{ss}(\chi)$ is open and $(G\times T)$-invariant in $Z$. Moreover, famous Mumford's construction
provides the categorical quotient $\pi:Z^{ss}(\chi)\to Z^{ss}(\chi)/\!/T$ of the $T$-action on $Z^{ss}(\chi)$, where
$$
Z^{ss}(\chi)/\!/T=\Proj(\bigoplus_{s\ge 0} \Of(Z)_{s\chi}).
$$
This shows that if a character $\chi$ defines a projective embedding $G/H\hookrightarrow X(\chi)$, then $X(\chi)$
is $G$-isomorphic to the quotient $Z^{ss}(\chi)/\!/T$.

Further, two characters $\chi_1$ and $\chi_2$ are called {GIT-equivalent}, if the corresponding sets 
$Z^{ss}(\chi_1)$ and $Z^{ss}(\chi_2)$ of semistable points coincide. 

\begin{proposition} \cite[Prop.~2.9]{BH}
Characters $\chi_1$ and $\chi_2$ are
GIT-equivalent if and only if $\sigma(\chi_1)=\sigma(\chi_2)$. 
\end{proposition}

In particular, the number of classes of GIT-equivalence is finite. The condition $Z^{ss}(\chi_1)=Z^{ss}(\chi_2)$
implies the canonical $G$-equivariant isomorphism of quotients: $Z^{ss}(\chi_1)/\!/T\cong Z^{ss}(\chi_2)/\!/T$.
Thus, we get

\begin{theo} \label{P2}
Assume that $H_1$ is a Grosshans subgroup of $G$. Then
\begin{enumerate}
\item a character $\chi$ defines
a projective embedding $G/H\hookrightarrow X(\chi) $ with small boundary if and only if $\chi\in C^{\circ}$;

\item if $\sigma(\chi_1)=\sigma(\chi_2)$, then the embeddings $X(\chi_1)$ and $X(\chi_2)$
are $G$-equivariantly isomorphic. 
\end{enumerate}
In particular, the number of isomorphism classes of projective $G/H$-embeddings with small boundary is finite.
\end{theo}

This theorem gives a partial answer to question (Q2). In the next section, we demonstrate that results of 
\cite{ArHa} give a complete answer to this question under some restrictions on the pair $(G,H)$.

\section{The generalized Cox's construction}
\label{sec:30}

 It is well known that any connected affine algebraic group $G$ is a semidirect product 
of a connected reductive subgroup $L$ and the unipotent radical $G^u$.
In turn, $L$ is an almost direct product $K\cdot S$ of a central torus $K$ and a semisimple subgroup $S$.
Consider a normal subgroup $\hat G:= SG^u$. Let $H$ be an epimorphic subgroup of $G$.

\begin{lemma} \label{L4}
The action $\hat G:G/H$ by left translations is transitive. 
\end{lemma}

\begin{proof}
We have to show that $\hat{G}H=G$. Consider the projection $\psi:G\to G/\hat G$ to the torus $G/\hat G$.
If $\psi(H)\ne G/\hat G$, then there exists a non-trivial character $\xi$ of the group $G$ whose
restriction to $H$ is trivial. Then $\xi$ may be considered as a non-constant regular
function on $G/H$. But $H$ is epimorphic, a contradiction. 
\end{proof}

 Further we assume that $G=\hat G$, or, equivalently, $\Chi(G)=0$. 
Lemma~\ref{L4} provides a partial compensation of this restriction. Moreover, we shall suppose that 
the class group $\Cl(G)$ is trivial. This may be achieved by replacing $G$ with its finite covering \cite[Prop.~4.6]{KKLV}.) 

\begin{definition}
We say that a subgroup $H\subset G$ is a {\it Grosshans extension}, if
$H$ is connected and $H_1=\cap_{\chi\in\Chi(H)}\Ker(\chi)$ is a Grosshans subgroup of $G$.
\end{definition}

If $H$ is a Grosshans extension, then $H/H_1$ coincides with the torus $T$. Under our conditions,
the algebra $\Of(G/H_1)$ is finitely generated and factorial. Indeed,
the condition $\Cl(G)=0$ is equivalent to factoriality of $\Of(G)$, connectedness of $H$ implies
connectedness of $H_1$ and the condition $\Chi(H_1)=0$ \cite[Prop.~3.13]{ArHa}, 
so the algebra $\Of(G)^{H_1}\cong\Of(G/H_1)$ is factorial \cite[Th.~3.17]{PV89}.

Following \cite{BH1}, we describe briefly a generalization of Cox's construction coming from toric geometry. 
This generalization delivers a realization of a wide class of algebraic varieties $X$ as a categorical quotient
for an action of so-called Neron-Severi torus $T$ on an open subset $\widehat X$ of an affine factorial variety $\overline{X}$. 
 
Let $X$ be an irreducible normal variety with a free finitely generated class group $\Cl(X)\cong\ZZ^k$.
Let us fix a subgroup $K$ in the group $\WDiv(X)$ of Weil divisors that projects isomorphically to $\Cl(X)$
under the natural projection $\WDiv(X)\to\Cl(X)$. Consider a graded sheaf of $\Of_X$-algebras on $X$:
$$
\mathcal{R}_X=\oplus_{D\in K} \Of(D), \ \ \text{где} \ \ \Of(D,U)=\{f\in\KK(X) : \div(f)+D\mid_U\ge 0\}. 
$$
The algebra of global sections of this sheaf  
$$
R(X)=\Gamma(\mathcal{R}_X,X)
$$
is factorial \cite{BH2}, \cite{EKW}, and is called the {\it total coordinate ring}, or the {\it Cox ring} of $X$.

Any homogeneous element $f\in R(X)$ of degree $D\in K$ defines a subvariety $Z(f)\subset X$, which is the support of the divisor
$\div(f)+D$. 

Suppose that the ring $R(X)$ is finitely generated. Then it corresponds to an affine factorial variety $\overline{X}:=\Spec(R(X))$.
Consider an open subset 
$$
\widehat{X}:=\bigcup_{f\in F(X)} \overline{X}_f\subset\overline{X},
$$
where $F(X)$ is a collection of homogeneous element of the ring $R(X)$ such that the open subset $X\setminus Z(f)$ is affine.
One may show that $\widehat{X}$ is isomorphic to the relative spectrum of the sheaf $\mathcal{R}_X$ over the variety $X$.
In particular, this relative spectrum is quasiaffine. 

The {\it Neron-Severi torus} of the variety $X$ is a torus $T$ whose lattice of characters is identified with the lattice $K$.
The $K$-gradings of the sheaf $\mathcal{R}_X$ and of the ring $R(X)$ define $T$-actions on the varieties $\widehat{X}$ and $\overline{X}$.
Clearly, the open embedding $\widehat{X}\subset\overline{X}$ is $T$-equivariant. The $T$-variety $\widehat{X}$ admits a categorical
quotient $\pi:\widehat{X}\to X$. The quotient morphism
$\pi:\widehat{X}\to X$ is called the {\it universal torsor} over the variety $X$. 

\medskip

In the case $X_0=G/H$ the Cox ring $R(X_0)$ coincides with the ring $\Of(G/H_1)$ \cite[Section~3]{ArHa}, and the finite generation condition is 
equivalent to the fact that $H_1$ is a Grosshans subgroup of $G$. The variety $\overline{X_0}$ is $Z=\Spec(\Of(G/H_1))$, and the universal torsor
is represented by a morphism $G/H_1\to G/H$, which is a quotient morphism for the right $T=H/H_1$-action on $G/H_1$. For arbitrary embedding
$G/H\hookrightarrow X$ with small boundary the variety $\overline{X}$ again coincides with $Z$, and $\widehat{X}$ is
an open $G$-invariant subset of $Z$ (in particular, it contains the open orbit $G/H_1$), and this open subset uniquely defines $X$.

Recall that a $G$-equivariant morphism
from an embedding $G/H\hookrightarrow X_1$ to an embedding $G/H\hookrightarrow X_2$ is a
$G$-equivariant morphism $\phi:X_1\to X_2$ that is identical on the open orbits.
Clearly, there exists at most one morphism between two given embeddings.
In \cite[Prop.~2.4]{ArHa}, we show that an equivariant morphism $\phi:X_1\to X_2$ between embeddings
corresponds to inclusion of open subsets $\widehat{X_1}\subseteq \widehat{X_2}$ of the variety $Z$. 
This leads to

\begin{theo}\cite[Th.~3.10]{ArHa} \label{T3}
Let $G$ be a connected algebraic group with $\Chi(G)=0$ and $\Cl(G)=0$,
and $H$ be an epimorphic subgroup of $G$, which is a Grosshans extension. Then 
\begin{enumerate}
\item two embeddings $G/H\hookrightarrow X(\chi_1)$ and $G/H\hookrightarrow X(\chi_2)$ are
$G$-isomorphic if and only if $\sigma(\chi_1)=\sigma(\chi_2)$;
\item there exists a $G$-equivariant morphism between embeddings $X(\chi_1)\to X(\chi_2)$ if and only if 
the cone $\sigma(\chi_2)$ is a face of the cone $\sigma(\chi_1)$.  
\end{enumerate}
\end{theo}

%%%%%%%%%%%%%%%%%%%%%%%%%%%%%%%%%%%%%%%%%%%%%%%%%%%%%%%%%%%%%%%%%%%%%%%%%%%%%%%%%%

\section{Geometry of embeddings with small boundary}
\label{sec:31}

 The methods of bunched rings theory \cite{BH1} allow to describe geometric properties
of normal varieties with a free finitely generated divisor class group and a finitely generated Cox ring
in combinatorial terms. This approach is based on the generalized Cox's construction discussed in the
previous section. Here we shall not present the general theory of bunched rings, but just formulate some
its results in conformity to projective embeddings with small boundary. Some information necessary for 
the translation of results may be found in \cite[Section~5]{ArHa}.

 Again, let $G$ be a connected algebraic group with $\Chi(G)=0$ and $\Cl(G)=0$,
and $H$ be an epimorphic subgroup of $G$, which is a Grosshans extension.
Since the GIT-fan $\Sigma(Z)$ is uniquely defined by the pair $(G,H)$, we shall denote it as $\Sigma(G/H)$.  

Following Section~\ref{sec:3}, we denote by $f_1,\dots,f_m$ a system of prime pairwise non-associate generators 
of the algebra $\Of(G/H_1)$ that are semiinvariants  with respect to the right $T=H/H_1$-action with weights
$\mu_1,\dots,\mu_m$. Denote by $E$ a lattice with a basis $e_1,\dots,e_m$ and by $Q:E\to\Chi(T)$ a projection 
sending $e_i$ to $\mu_i$. The same $Q$ will denote the corresponding projection of rational vector spaces
$E_{\QQ}\to \Chi(T)_{\QQ}$. Let $\gamma:=\cone(e_1,\dots,e_m)$ be a cone in $E_{\QQ}$ generated by
$e_1,\dots,e_m$. For any character $\chi\in\Chi(T)\cap\cone(\mu_1,\dots,\mu_m)$ denote by
$\cov(\chi)$ a collection of all faces $\gamma_0$ of the cone $\gamma$ such that $Q(\gamma_0)$ is the orbit cone
of some point $z\in Z=\Spec(\Of(G/H_1))$, $\chi\in Q(\gamma_0)^{\circ}$, and $\gamma_0$ is not a face
of any other face of $\gamma$ satisfying the same conditions. By $\lin(\gamma_0)$ denote the linear
span of a face $\gamma_0$ in the space $E_{\QQ}$.

By construction, the class group of the variety  $X(\chi)$ (coinciding with $\Cl(G/H)$) is identified with
the character lattice $\Chi(T)$ generated by $\mu_1,\dots,\mu_m$. The following proposition describes the Picard 
group $\Pic(X(\chi))$ as a sublattice in $\Chi(T)$. 

\begin{proposition} \cite[Prop.~7.1]{BH1}
$$
\Pic(X(\chi))=\bigcap_{\gamma_0\in\cov(\chi)} Q(\lin(\gamma_0)\cap E).
$$
\end{proposition}

Recall that a normal variety $X$ is called {\it locally factorial} if 
$\Cl(X)=\Pic(X)$, and {\it $\QQ$-factorial} if for any Weil divisor on $X$
some its multiple is a Cartier divisor. 

\begin{corollary} \label{corlf}
\begin{enumerate}
\item The variety $X(\chi)$ is locally factorial if and only if \linebreak
$Q(\gamma_0\cap E)$ generates the lattice $\Chi(T)$ for any $\gamma_0\in\cov(\chi)$. 
\item The variety $X(\chi)$ is $\QQ$-factorial if and only if the GIT-cone
$\sigma(\chi)$ has maximal dimension in $\Chi(T)_{\QQ}$.
\end{enumerate}
\end{corollary}

\begin{remark}
If the variety $\widehat{X(\chi)}$ is smooth, then the variety $X(\chi)$ is locally
factorial if and only if it is smooth, \cite[Prop.~5.6]{BH1}.
\end{remark}

 Denote by $\Eff(X)$, $\SAmple(X)$, and $\Ample(X)$ cones generated by divisor classes with an effective
representative, a base point free classes and ample classes respectively.

\begin{proposition} \cite[Prop.~7.2, Th.~7.3]{BH1}
$$
\Eff(X(\chi))=\cone(\mu_1,\dots,\mu_m), \ \ \SAmple(X(\chi))=\sigma(\chi), \ \ \Ample(X(\chi))=\sigma(\chi)^{\circ}.
$$
\end{proposition}

Finally, suppose that the ideal of relations between the elements $f_1,\dots,f_m$ in the ring $\Of(G/H_1)$
is generated by $K$-homogeneous polynomials $g_1,\dots,g_d$ with $d=m-\dim T-\dim X$. Then the results of
\cite[Section~8]{BH1} provide a formula for the canonical class $D_c$ of the variety $X(\chi)$: 

\begin{proposition} \label{prcc}
$$
D_c=\sum_{i=1}^d \deg(g_i)-\sum_{j=1}^m \mu_j.
$$
\end{proposition}

%%%%%%%%%%%%%%%%%%%%%%%%%%%%%%%%%%%%%%%%%%%%%%%%%%%%%%%%%%%%%%%%%%%%%%%%%%%%%%%

\section{Examples}
\label{sec:4}

In this section, $G$ is a connected simply connected semisimple algebraic group.

\begin{example}
Let us show that the fan $\Sigma(G/H)$ may have a complicated combinatorial structure.
Let $\chi_1,\dots,\chi_s$ be an arbitrary collection of non-zero elements of a lattice $M$
containing a basis of the lattice and generating a strictly convex cone $C$ in $M_{\QQ}$. 
Consider the set $\Omega$ of cones generated by all subsets of the set$\{\chi_1,\dots,\chi_s\}$. For any
character $\chi\in C\cap M$ define a cone
$$
\sigma(\chi):=\bigcap_{\omega\in\Omega, \chi\in\omega} \omega. 
$$
Proposition~\ref{P1} implies that the set of cones 
$$
\Sigma(\chi_1,\dots,\chi_s):=\{\sigma(\chi) : \chi\in C\cap M\}
$$
is a fan with support $C$. We claim that this fan is the GIT-fan $\Sigma(G/H)$ for some
homogeneous space $G/H$.

Consider a lattice $N$ with a basis 
$e_1,\dots,e_s$, and a surjective homomorphism $\phi:N\to M$ defined by
$e_i\to\chi_i$. If one identifies the lattices $M$ and $N$ with the character lattices of tori $T$ and $S$,
then the homomorphism $\phi$ represents $T$ as a subtorus of $S$. 

Let $S\subset B=SB^u$ be a maximal torus and a Borel subgroup of $G$.
We assume that the characters $e_1,\dots, e_s$ are identified with the fundamental weights
of the torus $S$ with respect to $B$. A subgroup $H:=TB^u$ is epimorphic in $G$, and the corresponding 
$H_1$ coincides with $B^u$. Using a description of the canonical embedding $G/B^u\hookrightarrow Z$ \cite[Th.5.4]{Gr}, 
one easily proves that any cone $\omega\in\Omega$ is an orbit cone for the $T$-variety $Z$
\cite[Prop.~4.4]{ArHa}. It follows that $\Sigma(G/H)=\Sigma(\chi_1,\dots,\chi_s)$.
\end{example}

\begin{example} \label{pr1}
Let $\gg$ be the tangent algebra of $G$, and $\nn\subset\gg$ be the null-cone of the adjoint module.
It is known that $\nn$ is irreducible, normal, contains a finite number of $G$-orbits, 
and all of them have even dimension \cite{Ko}. In particular, the complement of an open orbit $Ge\subset\nn$
has codimension $\ge 2$, and the connected component 
$G_e^0$ of the stabilizer $G_e$ is a (unipotent) Grosshans subgroup of $G$. Let
$H:=G_{\langle e\rangle}$ be the stabilizer of the line $\langle e\rangle$. 
Then $H^0$ is a Grosshans extension of the subgroup $H_1=G_e^0$ by a one-dimensional torus.
For example, if $G=\SL(4)$ then 
$$
H^0=\left\{ 
\begin{pmatrix}
t^3 & a & b & c \\
0 & t & a & b \\
0 & 0 & t^{-1} & a \\
0 & 0 & 0 & t^{-3} 
\end{pmatrix}
\ : \ 
t\in\KK^{\times}, \ a,b,c \in\KK \ 
\right\}.
$$
The projectivization $\PP(\nn)$ of the cone $\nn$ defines a projective embedding $G/H\hookrightarrow\PP(\nn)$ with small boundary.
Since the rank of the character lattice $\Chi(Н)$ equals one, this is a unique projective $G/H$-embedding with 
small boundary.
\end{example}

\begin{example}
If $G=\SL(2)$, then a Borel subgroup $B$ is a unique (up to conjugation) epimorphic subgroup,
the homogeneous space $G/B\cong\PP^1$ is projective and admits no non-trivial embeddings. 

Consider the case, where $G=\SL(3)$ and $H$ is connected. Here there are three projective homogeneous spaces: $G/B$,
$G/P_1\cong\PP(\KK^3)$, and $G/P_2\cong\PP((\KK^3)^*)$, where $P_1$ and $P_2$ are maximal parabolic subgroups.
Further, for the diagonal actions
$G:\PP(\KK^3)\times\PP(\KK^3)$ and $G:\PP((\KK^3)^*)\times\PP((\KK^3)^*)$ an open orbit has complement
of codimension $2$. This completes the list of projective embeddings $G/H\hookrightarrow X$ with small boundary, 
where $H$ has rank $2$.

Suppose that $H$ has rank $1$. Since $H$ is epimorphic, $\dim H\ge 3$ \cite{Bi}. 
Also $H$ contains no non-trivial semisimple subgroups, because otherwise one has
$\Chi(H)=0$, a contradiction with $H$ being epimorphic. Hence the subgroup $H$ is solvable.

If $H$ is regular, i.e., is normalized by a maximal torus of $G$, then there are three possibilities:

\smallskip

{\it Type 1.} \ 
$$
H=\left\{ 
\begin{pmatrix}
t^p & a & b \\
0 & t^q & c \\
0 & 0 & t^{-p-q}
\end{pmatrix}
\ : \ 
t\in\KK^{\times}, \ a,b,c \in\KK, \ p>0,\ p+q>0,\ (p,p+q)=1 \ 
\right\}.
$$
Here $H_1=B^u$ and 
$$
Z=\Spec(\Of(G/B^u))=\{(x_1,x_2,x_3,y_1,y_2,y_3) : x_1y_1+x_2y_2+x_3y_3=0\}\subset\KK^3\oplus(\KK^3)^*.
$$
We get a $2$-parameter family of projective embeddings with small boundary, namely
$X_{p,q}=(Z\setminus\{0\})/\!/\KK^{\times}$, where
$$
t\cdot(x_1,x_2,x_3,y_1,y_2,y_3)=(t^px_1,t^px_2,t^px_3,t^{p+q}y_1,t^{p+q}y_2,t^{p+q}y_3). 
$$
For each of them, the class group has rank one, and the Picard group is the subgroup of index $p(p+q)$. 

In this example, the variety $Z$ is a hypersurface. The weights of generators are $p,p,p,p+q,p+q,p+q$, 
and the weight of the relation is $2p+q$. Applying Proposition~\ref{prcc}, we find the canonical class of the variety $X_{p,q}$:
$$
D_c=2p+q-6p-3q=-4p-2q.
$$

\smallskip

{\it Type 2.} \ 
$$
H=\left\{ 
\begin{pmatrix}
t^p & 0 & b \\
0 & t^q & c \\
0 & 0 & t^{-p-q}
\end{pmatrix}
\ : \ 
t\in\KK^{\times}, \ b,c \in\KK, \ p,q>0 \ (p,q)=1 \ 
\right\}.
$$
Here we also have a 2-parameter family of embeddings corresponding to 
$$
X'_{p,q}=((\KK^3\oplus\KK^3)\setminus (0,0) /\!/\KK^{\times}
$$
with respect to the action
$$
t\cdot(x_{11},x_{12},x_{13},x_{21},x_{22},x_{23})=(t^px_{11},t^px_{12},t^px_{13},t^qx_{21},t^qx_{22},t^qx_{23}).
$$

\smallskip

{\it Type 3.} \ 
$$
H=\left\{ 
\begin{pmatrix}
t^p & a & b \\
0 & t^q & 0 \\
0 & 0 & t^{-p-q}
\end{pmatrix}
\ : \ 
t\in\KK^{\times}, \ a,b \in\KK, \ p,q>0,\ (p,q)=1 \ 
\right\}.
$$
Here the embeddings are obtained as in Type 2 via passage to the dual module $(\KK^3)^*$. 

Note that the spaces of embeddings for Types 2 and 3 are toric varieties.
One may write out their fans by a standard procedure.

\smallskip

Finally, there is one (up to conjugation) non-regular epimorphic subgroup:
$$
H=\left\{
\begin{pmatrix}
t & a & b \\
0 & 1 & a \\
0 & 0 & t^{-1}
\end{pmatrix}
\ : \ 
t\in\KK^{\times}, \ a,b \in\KK \ 
\right\},
$$
which also leads to a projective embedding $G/H\hookrightarrow X$ with small boundary. This embedding is a three-sheeted 
covering of the embedding of Example~\ref{pr1}.

\smallskip

Summing up, in the case $G=\SL(3)$ any homogeneous space $G/H$, where $H$ is connected epimorphic, admits exactly one projective embedding with
small boundary, and dimension of these embeddings varies from 3 to 5. 
\end{example}

\begin{example}
Let us give an example of an epimorphic subgroup $H$ of maximal rank in $G=\SL(4)$ such 
that $G/H$ admits many projective embeddings with small boundary.
Set
$$
H=\left\{ 
\begin{pmatrix}
t_1 & 0 & a & b \\
0 & t_2 & c & d \\
0 & 0 & t_3 & 0 \\
0 & 0 & 0 & t_4 
\end{pmatrix}
\ : \ 
t_1t_2t_3t_4=1, \ a,b,c,d \in\KK \ 
\right\}.
$$
Since the subgroup $H$ contains a maximal torus of $G$, $H$ is epimorphic if and only if
it is not contained in a proper reductive subgroup of $G$ \cite[Lemma 23.5]{Gr}. 
Note that $H^u$ coincides with the unipotent radical $P^u$ of a parabolic subgroup $P$.
It is known that in simple groups the unipotent radical of a parabolic subgroup is not contained in
a proper reductive subgroup \cite[Prop. 7]{Ar}. This proves that $H$ is epimorphic. On the other hand,
$H_1$ also coincides with $P^u$, thus it is a Grosshans subgroup
\cite[Th.16.4]{Gr}. 

In order to compute the fan $\Sigma(G/H)$, one needs some information on the variety $Z=\Spec(\Of(G/P^u))$.
Such varieties were studied in \cite{AT}. In our concrete case, \cite[Thm.~3.2]{AT}
provides the following explicit construction of the variety $Z$.
Let $e_1,e_2,e_3,e_4$ be the standard basis of $\KK^4$, and
$$
V_1=\Hom(\langle e_1,e_2\rangle, \KK^4), \ V_2=\Hom(\langle e_1\wedge e_2\rangle, \bigwedge^2\KK^4),
$$
$$ 
V_3=\Hom(\langle e_1\wedge e_2\wedge e_3, e_1\wedge e_2\wedge e_4\rangle, \bigwedge^3\KK^4) 
$$
be the spaces of linear maps regarded as $G$-modules with respect to the action on the image.
Consider $V=V_1\oplus V_2\oplus V_3$ as a $G$-module with the diagonal action,
and a point $x=(i_1,i_2,i_3)\in V$, where $i_j$ are the identical inclusions of the corresponding spaces. 
Then $Z=\overline{Gx}$. A $T$-action commuting with the $G$-action comes from the action on the argument
of a linear map. The weight semigroup of this action is generated by 
$$
\varepsilon_1, \ \varepsilon_2, \ \varepsilon_1+\varepsilon_2, \ \varepsilon_1+\varepsilon_2+\varepsilon_3, 
\ \varepsilon_1+\varepsilon_2+\varepsilon_4,
$$ 
where $\varepsilon_i((t_1,t_2,t_3,t_4))=t_i$. The relation $\varepsilon_1+\varepsilon_2+\varepsilon_3+\varepsilon_4=0$ implies
$\varepsilon_1+\varepsilon_2=(\varepsilon_1+\varepsilon_2+\varepsilon_3)+(\varepsilon_1+\varepsilon_2+\varepsilon_4)$.
The following diagram describes the location of the basic weights:

$$
\begin{array}{ccccc}
&& \varepsilon_1+\varepsilon_2+\varepsilon_3 && \\
&& | && \\
\varepsilon_1 & - & \varepsilon_1+\varepsilon_2 & - & \varepsilon_2 \\ 
&& | && \\
&& \varepsilon_1+\varepsilon_2+\varepsilon_4 &&
\end{array}
$$

\medskip
We claim that any cone generated by a subset of
$\{\varepsilon_1, \varepsilon_2, \varepsilon_1+\varepsilon_2+\varepsilon_3, \varepsilon_1+\varepsilon_2+\varepsilon_4\}$ is an orbit cone for the  $T$-action on $Z$. This claim may be checked directly.
For example, the cone generated by $\varepsilon_1, \varepsilon_2, \varepsilon_1+\varepsilon_2+\varepsilon_3$ is the orbit cone of the point
$$
y=((e_1\to e_1, e_2\to e_1), \ (e_1\wedge e_2\to 0), \ (e_1\wedge e_2\wedge e_3\to e_1\wedge e_2\wedge e_3, e_1\wedge e_2\wedge e_4\to 0)).
$$

Thus the fan $\Sigma(G/H)$ contains eighteen cones, and nine of them intersects the interior of $C$. By Theorem~\ref{T3}, the space $G/H$ admits
nine projective embeddings with small boundary, and the diagram of equivariant morphisms looks like:

$$ 
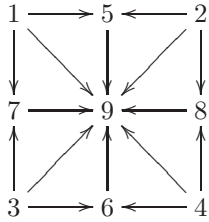
\xymatrix{
1 \ar[d] \ar[r] \ar[dr] & 5 \ar[d] & 2 \ar[l] \ar[d] \ar[dl]\\
7 \ar[r] & 9 & 8 \ar[l] \\
3 \ar[r] \ar[u] \ar[ur] & 6 \ar[u] & 4 \ar[l] \ar[ul] \ar[u]
}
$$ 

\end{example}

$\QQ$-factorial embeddings are represented by cones 1, 2, 3, 4. All of them are locally factorial, see Corollary~\ref{corlf}.

%%%%%%%%%%%%%%%%%%%%%%%%%%%%%%%%%%%%%%%%%%%%%%%%%%%%%%%%%%%%%%%%%%%%%%%%%%%%%%%%%

\end{document}